\RequirePackage{fix-cm}
\documentclass[a4paper,twoside,reqno,12pt]{amsart}

\usepackage{fixltx2e} 
\usepackage{a4wide,verbatim}
\usepackage[latin1]{inputenc}
\usepackage[english]{babel}
\usepackage{amsmath,amsfonts,amssymb,amsgen,amsbsy,dsfont,mathrsfs}
\usepackage[square,numbers]{natbib}
\usepackage{stmaryrd}

\usepackage{color}
\usepackage{graphicx}
\usepackage{float,here}
\usepackage{xcolor}
\usepackage{hyperref}
\usepackage{stmaryrd}

\usepackage{tikz}
\usetikzlibrary{positioning}
\usetikzlibrary{arrows}

\usepackage{mathtools}
\usepackage{multicol}

\hypersetup{
colorlinks=true, 
breaklinks=true, 
urlcolor= blue, 
linkcolor= blue, 
bookmarksopen=true, 
pdftitle={}, 
pdfauthor={Christian Bourdarias}, 
pdfsubject={} 
}

\newcommand{\half}{{\textstyle{1\over2}}}

\newcommand{\quat}{{\textstyle{1\over4}}}

\newcommand{\eqdef}{\stackrel{\text{\tiny{def}}}{=}}   

\newcommand{\cqfd}{{\nobreak\hfil\penalty50\hskip2em\hbox{}\nobreak\hfil
$\square$\qquad\parfillskip=0pt\finalhyphendemerits=0\par\medskip}}

\newtheorem{theorem}{Theorem}[section]

\newtheorem{proposition}{Proposition}[section]
\newtheorem{lemma}{Lemma}[section]
\newtheorem{definition}{Definition}[section]
\newtheorem{remark}{Remark}[section]
\newtheorem{exmp}{Example}[section]

\title[Regularizing effect for conservation laws]{\bf Regularizing effect for conservation laws \\ with a Lipschitz convex flux}

\author[Guelmame et al.]{Billel Guelmame, St\'ephane Junca and Didier Clamond}

\newcommand{\nfont}{\fontshape{n}\selectfont}
\address{({\nfont\textbf{Billel Guelmame}}) Universit\'e C\^ote d'Azur, Inria \& CNRS,   UMR 7351, LJAD,  
Parc Valrose, F-06108 Nice cedex 2, France.} 
\email{billel.guelmame@unice.fr}
\address{({\nfont\textbf{St\'ephane Junca}}) Universit\'e C\^ote d'Azur, Inria \& CNRS,  UMR 7351,  LJAD, 
Parc Valrose, F-06108 Nice cedex 2, France.} 
\email{stephane.junca@unice.fr}
\address{({\nfont\textbf{Didier Clamond}}) Universit\'e C\^ote d'Azur, CNRS, UMR 7351, LJAD, 
Parc Valrose, F-06108 Nice cedex 2, France.} 
\email{didier.clamond@unice.fr}

\date{\today}

\begin{document}

\maketitle

\begin{abstract}{This paper studies the smoothing effect for  entropy solutions 
of conservation laws with general nonlinear convex fluxes on $\mathds{R}$. Beside convexity, no additional 
regularity is assumed on the flux.  Thus, we generalize the well-known $\mathrm{BV}$ smoothing 
effect for $\mathrm{C}^2$ uniformly convex fluxes discovered independently by P. D. Lax \cite{Lax} and 
O. Oleinik \cite{O},  while in the present paper the flux is only locally Lipschitz. 
Therefore, the wave velocity can be dicontinuous and the one-sided Oleinik 
inequality is lost. This inequality is usually the fundamental tool to get a sharp regularizing 
effect for the entropy solution. 
We modify the wave velocity in order to get an Oleinik inequality useful for the wave front 
tracking algorithm. Then, we prove that the unique entropy solution belongs to a generalized $\mathrm{BV}$ 
space, $\mathrm{BV}^\Phi$.
} 
\end{abstract}

\medskip


\medskip



\tableofcontents

\section{Introduction}

This paper is about the regularization effect on the unique entropy solution of the scalar hyperbolic 
conservation law 
\begin{equation}\label{cl}
u_t\ +\ f(u)_x\ =\ 0, \qquad u(0,x)\ =\ u_0(x), \qquad  M=\|u_0\|_\infty.
\end{equation}
The initial datum $u_0$ belongs to $\mathrm{L}^\infty(\mathds{R},\mathds{R})$. 
In (\ref{cl}), $f$ is a nonlinear convex flux on the whole real line, thence $f$ is Lipschitz on 
$[-M,M]$. The regularity of $u$ for positive time $t$ depends on the nonlinearity of $f$ on $[-M,M]$. 
(For a linear flux, the solution is nothing but a translation of the initial datum  with a constant speed, 
so no regularity is enforced by the equation (\ref{cl}).) To obtain a smoothing effect, the 
following Tartar condition \cite{Tartar} is needed:
\begin{equation}\label{Tartar}
\mbox{\em There\ are \ no\ non-trivial\ interval\ where\ f\/\ is\ affine}.
\end{equation}
Here, the flux being nonlinear and convex on $\mathds{R}$, it is strictly convex and thus it 
necessarily satisfies the condition \eqref{Tartar}.

In \cite{Lax,O}, both Lax and Oleinik prove that for an uniform convex flux $f$ such that  
$f''>c>0$ for some constant $c$ (as, e.g., for the Burgers equation), the solution $u(t,\cdot)$ 
is in $\mathrm{BV_{loc}}$, for all time $t>0$. (Definitions of the various $\mathrm{BV}$ 
spaces, spaces of functions of bounded variation, can be found below and also in \cite{AFP,MO}.) 
This result is no longer true for flatter fluxes \cite{Cheng},  such as $f(u)=|u|^3$ and $f(u)=u^4$. 
The solution regularities in SBV and Sobolev spaces are obtained in \cite{AGV,Jabin}.
To obtain more information on the regularity of $u$, generalized BV spaces, $\mathrm{BV}^s$ and 
$\mathrm{BV}^\Phi$, are needed. The regularity in those spaces implies the right 
regularity in Sobolev spaces, as well as the left and right traces for shock waves. For smooth 
fluxes with a polynomial degeneracy (e.g., $f(u)=|u|^3$, $f(u)=u^4$), the solution $u(t,\cdot)$ 
belongs to $\mathrm{BV}^s_\mathrm{loc}$ in space for $t>0$ (see the last paragraph of Section \ref{sec2} and 
\cite{BGJ,CJJ}, and see \cite{BM,Elio} for non-convex fluxes). This kind of regularity is still true 
for a  $\mathrm{C}^1$ convex flux in a bigger generalized BV space, $u(t,\cdot) \in 
\mathrm{BV}^\Phi_\mathrm{loc}$, $t>0$ with a convex function $\Phi$ depending on the nonlinearity 
of $f$ \cite{CJLO,Elio}. 

In this paper, we show that this last result \cite{CJLO} remains true for all convex fluxes 
$f$ on $\mathds{R}$ satisfying the condition \eqref{Tartar}, without requiring 
$f$ to be in $\mathrm{C}^1$.
Such  a flux can appear in applications, such as in traffic flow model \cite{W} 
with a concave flux. 
If $f$ is a strictly convex flux, then the (necessarily increasing) velocity 
\begin{equation}\label{velocity}
a(u)\ \eqdef\ f'(u)
\end{equation} 
exists 
almost everywhere. The set of discontinuity of $a$ is countable, the left and right limits 
$a^-(u) \leqslant a^+(u)$  existing everywhere. Thanks to the maximum principle, the entropy 
solution $u$ takes values only in $[-M,M]$, hence $a$ is bounded on $[-M,M]$. In the case of 
$\mathrm{C}^1$ convex fluxes, the simplest proof (see \cite{CJJ} after \cite{O}) is based on 
the fundamental one-sided Oleinik inequality \cite{Hoff}, 
\begin{equation}\label{Ol}
a(u(t,x))\,-\,a(u(t,y))\ \leqslant\ (x-y)/t \quad \mathrm{a.e.}\quad x\,>\,y,
\end{equation}
which implies that $a(u)$ is a BV function and then $u$ belongs to a $\mathrm{BV}^\Phi$. 
Unfortunately, this inequality is no longer true for convex Lipschitz fluxes. Indeed, first, 
$a(u)$ is not well defined because $a$ is not continuous and, second, the Oleinik inequality 
 is not true almost everywhere, as shown in Example \ref{ex:Oleinik-wrong} below. 
To our knowledge, the loss of the Oleinik inequality appears in the classical literature 
of conservation laws only in \cite{Hoff}, for a piecewise linear flux. Note that, though not always 
true, the Oleinik inequality is true on a large subset of $\mathds{R^+}\times\mathds{R}$. We prove 
in this paper that this is enough to still obtain the smoothing effect in the right $\mathrm{BV}^\Phi$ 
space with a modified wave velocity and a wave front tracking algorithm for scalar conservation 
laws \cite{Dafermos72}.

The paper is organized as follows. 
In Section \ref{sec2}, the function $\Phi$ is built to state the main Theorem.
The loss of the Oleinik inequality and the resulting difficulty to prove the main Theorem \ref{maintheorem} is discussed in Section \ref{sec2bis}. The Section \ref{sec3} studies the approximate 
Riemann problem and a modified Oleinik inequality. 
The Section \ref{sec4} is devoted to obtaining a BV estimates on the modified velocity 
by the wave front tracking algorithm. The main result is proved in Section \ref{sec5}.

\section{The main result}\label{sec2}

In this section, definitions of weak entropy solutions and  $\mathrm{BV}^\Phi$ spaces are  recalled   
and the function $\Phi$ related to the smoothing effect is built. Then, the smoothing effect is stated 
in Theorem \ref{maintheorem}.

\begin{definition}
$u$ is called a weak solution of \eqref{cl}, if for all  smooth functions 
$\theta$ with a compact support, i.e., for  $\theta\in\mathcal{D}(\mathds{R}^+\!\times\mathds{R})$
\begin{equation}\label{ws}
\int_\mathds{R}\int_{\mathds{R}^+}\left[\,u(t,x)\,\theta_t(t,x)\,+\,f(u(t,x))\,\theta_x(t,x)\,\right]
\mathrm{d}t\,\mathrm{d}x\ +\ \int_\mathds{R} u_0(x)\,\theta(0,x)\,\mathrm{d}x\ =\ 0.
\end{equation}
\end{definition}
For a given $u_0\in\mathrm{L}^\infty$, the equation \eqref{ws} has at least one weak solution 
\cite{Bressan, DafermosBook}, the uniqueness being ensured by the Kruzkov entropy conditions:
\begin{definition}{\bf[Kruzkov entropy solutions]}
A weak solution of \eqref{ws} is called an entropy solution if for all positive $\theta\in
\mathcal{D}(\mathds{R}^{+*}\!\times\mathds{R})$ and for all convex functions $\eta\in\mathrm{C}^1$ 
and with $F\eqdef\int f'(u)\eta'(u)\,\mathrm{d}\/u$ (primes denoting the derivatives), the following 
inequality holds:
\begin{equation}\label{entropy}
\int_\mathds{R}\int_{\mathds{R}^+}\left[\,\eta(u(t,x))\,\theta_t(t,x)\,+\,F(u(t,x))\,\theta_x(t,x)\,
\right]\mathrm{d}t\,\mathrm{d}x\ \geqslant\ 0.
\end{equation}
In addition, $u$ has to belong to  $\mathrm{C}^0([0,+\infty[,\mathrm{L^1_{loc}}(\mathds{R}))$.
\end{definition}
The functional space $\mathrm{BV}^\Phi$ \cite{MO} is defined as follow. 
\begin{definition}{\bf [$BV^\Phi$ spaces]}\label{TVPhi}
Let $\Phi$ be a convex function such that $\Phi(0)=0$ and $\Phi(h)>0$ for $h>0$, the 
total $\Phi$-variation of $v$ on $K \subset \mathbb{R}$ is
\begin{equation}
\mathrm{TV}^\Phi\/v\ \{K\} =\ \sup_{p \in \mathcal{P}}\ \sum_{i=2}^n \Phi(|v(x_i)-v(x_{i-1})|)
\end{equation}
where $\mathcal{P}=\{ \{x_1, \cdots ,x_n \},\ x_1<\cdots<x_n  \} $ is the set of all subdivisions  
of\/ $K$. The space $\mathrm{BV}^\Phi$ is defined by\/  
$\mathrm{BV}^\Phi=\{v, \, \exists \lambda >0, \,  \mathrm{TV}^\Phi (\lambda v)<  \infty \}$.
\end{definition}
The $\mathrm{BV}^\Phi$ space, a generalization of the $\mathrm{BV}$ space, is the space of functions 
with generalized bounded variations. Our goal here is to try to  construct the  optimal convex function $\Phi$, such that $u$ 
is in $\mathrm{BV}^\Phi_\mathrm{loc}$. 
 That means choosing $\Phi$ to obtain the  optimal  space $\mathrm{BV}^\Phi$ in order to characterize the regularity of the entropy solution. 
 A discussion on the optimal choice is done below in Remark \ref{optimality}.

The one-sided Oleinik inequality is directly linked with the increasing variation of an entropy solution. 
Such variation is defined  with $y^+ = \max(y,0)$ as follow with the same notations as in the previous 
definition \ref{TVPhi}. 
\begin{align}\label{TVPhi+}
\mathrm{TV}^{\Phi+}\/v\ \{K\} =\ \sup_{p \in \mathcal{P}}\ \sum_{i=2}^n \Phi((v(x_i)-v(x_{i-1}))^+)
\end{align}

In order to build the function $\Phi$, a definition of the generalized inverse of non decreasing functions is needed
\begin{definition}{\bf[Generalized inverse]}\label{GI}
Let $g$ be a non decreasing function from $\mathbb{R}$ to $\mathbb{R}$, the generalized inverse of $g$ 
is defined on $ g(\left [-M,  M \right])$ as following 
\begin{equation}\label{inverse}
g^{-1}(y)\ \eqdef\ \mathrm{inf} \{x \in \mathds{R}, \, y \leqslant g(x)\},
\end{equation}
\end{definition}

\begin{remark}\label{leftinverse} It's obvious from the definition that 
\begin{equation}
(g^{-1} \circ g)(x) \leqslant x, \; \forall x.
\end{equation}

\end{remark}
The usual properties of the generalized inverse can be found in \cite{inverse}.
\begin{proposition}[see Proposition 2.3 in \cite{inverse}]\label{properties} {\color{white} Billel }
\begin{itemize}
\item[1-] If $g$ is continuous then $g^{-1}$ is a strictly increasing function and
\begin{equation}
(g \circ g^{-1})(y)=y,\; \forall y
\end{equation}
\item[2-] If $g$ is strictly increasing then $g^{-1}$ is continuous and 
\begin{equation}
(g^{-1} \circ g)(x) = x,  \; \forall x,y
\end{equation}
\end{itemize}
\end{proposition}

Now, the function $\Phi$ is built. Notice that the velocity $a$ defined in \eqref{velocity} is strictly increasing due to the strict convexity of the flux $f$, then Proposition \ref{properties} ensure that the generalized inverse (see Definition \ref{GI}) of the velocity
\begin{equation}
b =\ a^{-1},
\end{equation}
is continuous and non decreasing. Note  also that the function $b$ is constant on $[a^-(u),a^+(u)]$ when $a$ is discontinuous 
at $u$.

Let $\omega[b]$ be the modulus of continuity of $b$, i.e., 
\begin{equation}\label{modulus}
\omega[b](h)\ \eqdef\ \sup_{\substack{|x-y|\leqslant h \\ x,y \in a([-M,M])}}|b(x)-b(y)|,
\end{equation}
and let $\phi$ be the generalized inverse of $\omega[b]$, i.e., $\phi(y)=\inf \{x \in \mathbb{R}, 
\, y \leqslant \omega[b](x)\}$. We denote $\Phi$ the lower convex envelope of $\phi$, that is related 
to the nonlinearity of the flux via the velocity. Let us write, in a concise way, the definition 
of $\Phi$. This function is the key ingredient to define a suitable functional space describing the 
regularity of entropy solutions. 
\begin{definition}{\bf [Choice of $\Phi$]}\label{defPhi}
$\Phi$ is the lower convex envelope of the generalized inverse of the modulus of 
continuity of the generalized inverse of the velocity, 
\begin{equation}\label{def:Phi}
\Phi\ \eqdef\ \mathrm{lower\ convex\ envelope\ of}\left.  \left( \omega[a^{-1}] \right)^{-1} \right. .
\end{equation} 
\end{definition}
\begin{remark}{\bf[On the optimality of the convex function $\Phi$]}\label{optimality}
This definition generalizes the ones given in \cite{BGJ,CJLO} for a discontinuous velocity.  
This is the optimal choice for a flux with a power law degeneracy, as proved in  e.g. \cite{CJ,GJ}, for the convex power flux $f(u)= |u|^{1+p}/(1+p)$, $p \geqslant 1$, $\Phi(u)= |u|^p = |a(u)|$. 
Indeed, when the velocity is convex for $u>0$ and is an odd function, then $\Phi(u)=  |a(u)|$ \cite{CJLO}.
\end{remark}

$\mathrm{C}_{w}$ denoting the space of continuous functions with modulus of continuity $w$, 
the goal of the paper is to prove the Theorem: 
\begin{theorem}{\bf[Regularizing effect in $\mathrm{BV}^\Phi$]}\label{maintheorem}
Let $f$ be a strictly convex flux on\/ $\mathds{R}$, $u_0 \in \mathrm{L}^\infty$ and $u$ being the 
unique entropy solution of \eqref{cl}, then $u(t,\cdot) \in \mathrm{BV}^\Phi_\mathrm{loc}$, i.e., 
for all $[\alpha, \beta]\subset\mathds{R}$,
\begin{align}\label{mt+}
\mathrm{TV}^{\Phi+} u(t,\cdot) \{[\alpha,\beta] \}\ &\leqslant\, \left( \beta-\alpha \right) t^{-1}, \\
\mathrm{TV}^{\Phi} u(t,\cdot) \{[\alpha,\beta] \}\ &\leqslant\  2 \left( \|a(u_0)\|_\infty  + 
(\beta-\alpha)\,t^{-1} \right),
\end{align}
Moreover, if $u_0$ is compactly supported, then there exists $C>0$ such that
\begin{equation}\label{mt}
\mathrm{TV}^\Phi u(t,\cdot) \{\mathds{R} \}\ \leqslant\ C\left(1+t^{-1}\right), \qquad 
\end{equation}
and, in addition, for all $\tau>0$ 
\begin{equation}\label{reg-time}
u \in   \mathrm{C}_{w_\tau} \left(]\tau,+\infty[,
\mathrm{L}^1_\mathrm{loc}\right) \qquad \mathrm{where} \qquad 
w_\tau(y)=\Phi^{-1}(C\left(1+\tau^{-1}\right) y ).
\end{equation}
\end{theorem}
Inequality \eqref{mt+} is the natural way to recover the one-sided Oleinik inequality. 
When $\Phi$ is the identity function, so $\mathrm{TV}^\Phi u= \mathrm{TV} u$,  the inequality 
\eqref{mt} is the classical one for uniformly convex smooth flux. 
The regularity in time \eqref{reg-time} is proven in Section \ref{sec5}.
\\

Theorem  \ref{maintheorem} covers all previous results (with a different proof) on the 
smoothing effect for a strictly convex flux, the $\mathrm{C}^2$ case being considered 
in \cite{BGJ,CJJ,Pierre,CJ,Jabin, Lax,Elio,O} and the $\mathrm{C}^1$ case being treated 
in \cite{CJLO}. All these proofs make use, directly or indirectly, of the Oleinik 
inequality (\ref{Ol}).  The proof of Theorem \ref{maintheorem} here is necessarily 
more complicated due to the loss of the Oleinik inequality. This crucial point is 
discussed in details  in the next section.

\section{Notes on the Oleinik inequality for discontinuous wave speeds}
\label{sec2bis}

The following example shows that the  Oleinik inequality \eqref{Ol} is no longer true everywhere 
when the velocity $a$ defined in \eqref{velocity} is not continuous on $[-M,M]$. The Oleinik inequality requires the velocity is defined everywhere. For this purpose, the velocity can be defined everywhere as the mean 
of its left and right limits with a weight $\lambda \in [0,1]$, i.e.,
\[
\bar{a}_\lambda(x)\ \eqdef\ \lambda\, a^{+}(x)\ +\ (1-\lambda)\,a^{-}(x).
\] 
Now, a key result about the Oleinik inequality for the proof of Theorem is stated with this 
mean velocity. 

Let us consider the Riemann problem, that is a Cauchy problem with a piecewise constant datum  
$u_0(x)= u_l$ for $ x < 0$ and  $u_0(x)= u_r$ for $ x > 0$. 
The  Oleinik inequality is clearly true for a shock wave, but it is not always valid 
for a rarefaction wave. 
\begin{proposition}{\bf [One-sided Oleinik inequality]} \label{prop:whereOlga}
For a Riemann problem producing a rarefaction wave --- i.e., $u_l < u_r$ --- if $x/t, y/t \in 
\bar{a}_\lambda([-M,M])$ then the Oleinik inequality holds
\begin{equation}\label{Olbar}
\bar{a}_\lambda(u(t,x))\ -\ \bar{a}_\lambda(u(t,y))\ \leqslant\ (x-y)/t \qquad \mathrm{a.e.}\quad x\,>\,y. 
\end{equation} 
\end{proposition}
The set $\bar{a}_\lambda([-M,M])$ is not an interval since $a$ is not continuous. This is a reason for the 
loss of the Oleinik inequality. Moreover, the solution is constant where $a$ is not defined, as shown 
in Example \ref{ex:Oleinik-wrong} below. When the velocity 
is continuous this problem disappears, as proved in the next section.

\begin{remark}
The Oleinik inequality is true a.e. with a velocity chosen as 
\begin{equation}\label{Oltrue}
a^-(u(t,x))\ -\ a^+(u(t,y))\ \leqslant\ (x-y)/t \qquad \mathrm{a.e.}\quad x\,>\,y.
\end{equation}
But, it is less useful to get the\/ $\mathrm{BV}$ estimate for this velocity. This is the key 
point to prove the\/ $\mathrm{BV}^\Phi$ regularity of the entropy solution.
It should also be noticed that the Oleinik inequality \eqref{Oltrue} can be invalid if the exponent signs of $a(.)$
are exchanged. 
\end{remark}
From now on, we denote $\xi=x/t$ and $\eta=y/t$ for brevity.
\begin{exmp} \label{ex:Oleinik-wrong}
{\bf[The one-sided Oleinik inequality is not always valid]}
Consider 
\begin{align}
f(u)\, =\,  u^2 + |u|,
%
\qquad
a(u)\,=\,
2u + \mathrm{sign}(u)
\end{align}
and $u_0(x)= \mathrm{sign}(x)$. The entropy solution of \eqref{cl} is $u(t,x)=U(\xi)$ with
\vspace{5mm}

\begin{tabular}{ccc}
\(
\hspace{-10mm} U(\xi)\, =\, 
\left\{
\begin{array}{cl}
-1            & \xi \leqslant -3, \\
\half(\xi+1) \quad & -3 \leqslant \xi \leqslant -1, \\
0             & -1 \leqslant \xi \leqslant 1, \\
\half(\xi-1)  & 1 \leqslant \xi \leqslant 3, \\
1             & 3 \leqslant \xi.
\end{array}
\right.
\)
&
\begin{tikzpicture}[thick, transform canvas={scale=0.8}, shift={(1.8,0)}]

\draw[thick,->] (-2.5,0) -- (2.5,0) node[anchor=north east] {$\xi$};
\draw[thick,->] (0,-2) -- (0,2) node[anchor= east] {$U$};

\draw[line width= 0.5 mm, blue ] (-2,-1) -- (-1.5,-1);
\draw[line width= 0.5 mm, blue ] (-1,0) -- (-1.5,-1);
\draw[line width= 0.5 mm, blue ] (-1,0) -- (1,0);
\draw[line width= 0.5 mm, blue ] (2,1) -- (1.5,1);
\draw[line width= 0.5 mm, blue ] (1,0) -- (1.5,1);

];

\end{tikzpicture}
&
\begin{tikzpicture}[thick, transform canvas={scale=0.8}, shift={(1,0)}]
\draw[thick,->] (4,-1.5) -- (10,-1.5) node[anchor=north west] {$x$};
\draw[thick,->][color=red] (7,-0.25-1.5) -- (7,3.5-1.5) node[anchor= east] {$t$};
\draw[thick] (4,0-1.5) -- (9,0-1.5) node[anchor=south] {$u=1$};
\draw[thick] (10,0-1.5) -- (5,0-1.5) node[anchor=south] {$u=-1$};
\draw[thick] (7,2-1.5) node[anchor= north] {\Large $u=0$};
\draw[thick][] (9.2,0.6-1.5) node[anchor=south east, transform canvas={scale=0.8}, 
shift={(3,0)}]{$u=\frac{\xi-1}{2} $};
\draw[thick][] (-2.2,0.6-1.5) node[anchor=south west, transform canvas={scale=0.8}, 
shift={(8,0)}]{$u=\frac{\xi+1}{2} $};
\draw[thick][color=blue] (7,0-1.5) -- (7-3,1.1-1.5) node[anchor=south]{$\xi=-3$};
\draw[thick][color=blue] (7,0-1.5) -- (7-9/4,3-1.5 ) node[anchor= west]{$\xi=-1$};
\draw[thick][color=blue] (7,0-1.5) -- (7+9/4,3-1.5 ) node[anchor= east]{$\xi=1$};
\draw[thick][color=blue] (7,0-1.5) -- (10,1.1-1.5) node[anchor=south]{$\xi=3$};
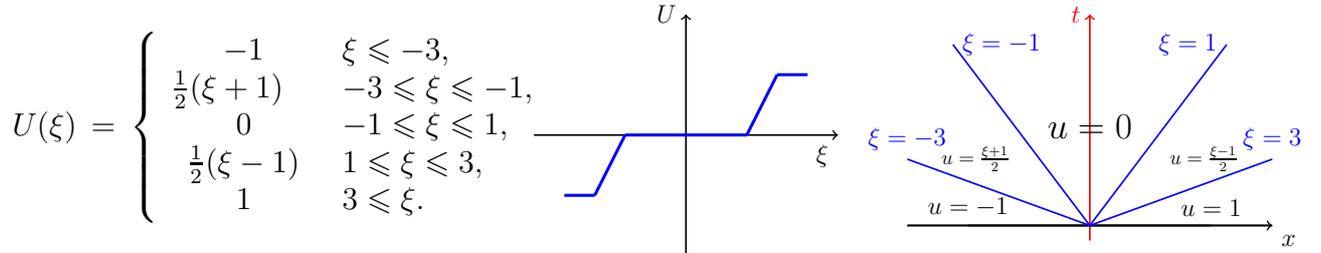
\node at (0,-2.2) {\textbf{Fig. 1.} The solution of the problem.};
\end{tikzpicture}
\end{tabular}

\begin{center}
\end{center}
\vspace{8mm}
Considering $t>0$, the Oleinik inequality is not satisfied and 
 $\bar{a}_\lambda(u(t,x))-\bar{a}_\lambda(u(t,y))>(x-y)/t$ 
in the following cases
\begin{itemize} 
\item if\/ $\lambda = 0$ and $-1<\eta<1<\xi<3$; 
\item if\/ $\lambda = 1$ and  $-3<\eta<-1<\xi<1$;
\item if\/ $\lambda \in ]0,1[ $ and  $\left(\ 2\lambda-1<\eta<1<\xi<3\ \mathrm{or}\  
-3<\eta<-1<\xi<1-2\lambda\ \right)$.
\end{itemize} 
\end{exmp}

\begin{remark} 
If we change the flux $f$ in Example \ref{ex:Oleinik-wrong} by $f(u)=u^2+u+|u|$ and if 
$\alpha<0<\beta$, then $\bar{a}_\lambda(u(t,\beta))-\bar{a}_\lambda(u(t,\alpha)) \geqslant 1$, where the 
inequality \eqref{mt+} remains true. 
\end{remark}
\begin{remark}  The converse of Proposition \ref{prop:whereOlga} is false in general.
For instance, it is sufficient to take $x,y$ in $]-t,0[$ or in $]0,t[$ in Example 
\ref{ex:Oleinik-wrong}.
\end{remark}

\begin{remark}
The function $b$ is the generalized inverse of 
 $\bar{a}_\lambda$      for all $\lambda$.  That is to say that $b$ does not depend on $\lambda$:  $$ \forall \lambda \in [0,1], \;   \forall u, \qquad  b(\bar{a}_\lambda(u))=u.$$
\end{remark}
 Indeed, the main result (Theorem \ref{maintheorem}) does not depend on the choice of $\lambda$. 
Thus, we can take $\lambda=1/2$ to fix the notation in all the sequel: 
\begin{equation}
\bar{a}(x)\ \eqdef\ \bar{a}_{1/2}(x)\ =\ \half\,[\,a^+(x)\,+\,a^-(x)\,].
\end{equation}


At this stage, it is important to outline the main difficulties for proving Theorem \ref{maintheorem}.
 These difficulties result from the discontinuity of the velocity $a(u)$
(yielding the loss of the Oleinik inequality).  For instance, consider the case when 
$a(u)$ is discontinuous at $u=u^\#$ with the jump  $$\llbracket a \rrbracket (u^\#)= a^+(u^\#)
-a^-(u^\#)>0.$$ Let $u(x,t)$ be a rarefaction wave, non-decreasing 
solution of the Riemann problem with initial data $u_0(x)=u_l$ for $x<0$ and $u_0(x)=u_r$ 
for $x>0$ with  $u_l < u^\# < u_r $.  The solution $u(x,t)$ is flat, with value $u=u^\#$  
on the interval $x\in[a^-(u^\#)\times t;a^+(u^\#)\times t]$, see figure 1. The size of this interval,  
at time $t$, is exactly 
\[
\Delta x\ =\ t\times \llbracket a \rrbracket (u^\#).
\] 
At first sight, it seems a good case where the Oleinik inequality is actually an equality. 
However, it is not the case for the two following reasons.

First, $u$ being constant on the flat part (with $u \equiv u^\#$), there are of course no 
variations of $u$ on this part, while there is a variation of $a$ equals to the jump of $a$ 
at $u^\#$. This shows that the variations of $a$ are bad indicators of the variations of $u$. 
Usually, the total variation of $u$  is controlled by the total variation of $a(u)$ \cite{CJJ}. 

Second, since $u$ is constant on this part, 
the shock wave penetrates the flat part, reducing the size of this part, i.e. $$ \Delta x < 
t\times \llbracket a \rrbracket (u^\#).\footnote{Note that it is the converse of the Oleinik inequality.}
$$ In other 
words, the jump of $a(u)$ at $u^\#$ does not represent well the 
size of the flat part, which is problematic as already mentioned. As a consequence, the total 
variation of $a(u)$ is not controlled in the present work. It is an important difference with 
the case of smooth fluxes, where $a(u)$ is known to be in $\mathrm{BV}$ \cite{Cheng2},  at least $\text{C}^2$, with precise assumptions and 
counter-examples given in \cite{Elio}. Here, assuming only that the convex flux is Lipchitz, it is not clear whether $a(u)$ belongs to $\mathrm{BV}$ or not. 
\medskip \\

If the velocity has only one discontinuity, it is easy to overcome this difficulty.  Consider 
$$ a(u)=\tilde{a}(u) +\llbracket a(u^\#)\rrbracket  H(u- u^\#), $$ where $H$ is the Heaviside function; that is 
to say $ \tilde{a}$ is the continuous part of the velocity $a(u)$. The generalised inverses 
of $a$ and $\tilde{a}$ having the same modulus of continuity, they define the same space 
$\mathrm{BV}^\Phi$. Moreover, $\tilde{a}$ is easy to estimate in $\mathrm{BV}$.  Thus, the 
$\mathrm{BV}^\Phi$ regularity of $u$ follows as in \cite{CJJ}.  This simple case shows again 
that the variation of $a$ through its discontinuity is useless to capture the regularity of $u$. 
Moreover, for Example \ref{ex:Oleinik-wrong} above,  the regularity of $u$ is simply 
$\mathrm{BV}$ since  $\tilde{a}$ corresponds to a Burgers flux.
\\

Removing the discontinuity of 
$a$ can be done only if $a$ has finitely many points of discontinuity. This is not possible 
in general. Consider the example with the following velocity 
$$ a(u) = \sum_n  2^{-n}  H( u-r_n), $$ 
where the sequence $(r_n)$ takes all the values of rational numbers. This velocity corresponds 
to a strictly convex flux where the second derivatives has only an atomic part. Thus, $u$ is 
regularised in some $\mathrm{BV}^\Phi$. Removing all the discontinuities of $a$ is not a good 
idea here since the corresponding flux is flat and does not correspond to any smoothing effect. 
\\
 
The way we solve this difficulty is to keep the Oleinik inequality by introducing a new velocity, 
called $\chi$ below. The total variation of $\chi$ can be estimate geometrically by the mean 
of characteristics through a wave front tracking algorithm. Moreover, the variation of $\chi$ 
corresponds exactly to the generalised variation of the entropy solution $u$.
\\

In order to prove Theorem \ref{maintheorem}, the wave front tracking algorithm will be used.
To do so, $u_0$ is approximated by a sequence of step functions and thus the Riemann problem for 
each sequence can be solved \cite[Lemma 3.1]{Dafermos72}. A Riemann problem with a discontinuous 
velocity is expounded in the next section. 

\section{Riemann problem}\label{sec3}

This Section is devoted to the Riemann problem. The shock wave is solved as usual, but solving the 
rarefaction wave is more complicated. For this purpose, the flux $f$ is approached by piecewise quadratic 
$\mathcal{C}^1$ fluxes, in order to show that, for a rarefaction wave, the solution is given by 
$b\left( x/t \right)$; this is the classic formula for smooth fluxes where $b$ is the inverse of the velocity. 
We extend this formula when $b$ is the generalized inverse of a discontinuous velocity. The second part 
deals with a piecewise linear flux $f_\varepsilon$ \cite{Dafermos72}, which gives a modified Oleinik 
inequality up to  a small error.

\subsection{The exact solution of the Riemann problem}\label{exactsol}

The Riemann problem consists in solving $ u_t+f(u)_x=0$ with the initial condition
\begin{equation}\label{Rm} 
u(0,x)\ =\ u_l\quad\mathrm{for}\quad x<0 \qquad\mathrm{and}\qquad 
u(0,x)\ =\ u_r\quad\mathrm{for}\quad x>0. 
\end{equation}
If $u_l>u_r$, the solution generates a shock with a speed given by the Rankine-Hugoniot relation 
$s=\llbracket f(u) \rrbracket/\llbracket u \rrbracket$, where $\llbracket u \rrbracket=u_r-u_l$ 
is the jump of $u$. The entropy solution is 
\[ 
u(t,x)\ =\ u_l\quad\mathrm{for}\quad x<s\/t \qquad\mathrm{and}\qquad 
u(t,x)\ =\ u_r\quad\mathrm{for}\quad x>s\/t. 
\]
The interesting case is $u_l<u_r$ because the Oleinik inequality is not always true in this case. 
The solution has a non decreasing rarefaction wave between $x=a^+(u_l)t$ and $x=a^-(u_r)t$, given 
by the following proposition:
\begin{proposition}\label{prop:rarefaction}
Let $u$ be the entropy solution of the Riemann problem with $u_l<u_r$. For $\xi=x/t$ the solution is
\begin{equation} \label{b-rarefaction}
u(t,x)=\begin{cases} 
      u_l & \xi < a^-(u_l) , \\
      b\left(\xi \right) & a^-(u_l)<\xi<a^+(u_r) , \\
      u_r & \xi > a^+(u_r) . 
   \end{cases} 
\end{equation}
\end{proposition}

\begin{remark}
A similar formula  for systems with Lipschitz fluxes is given in \cite[Th. 3.3, p. 279]{LeFloch}.  
Another formula is proposed by Bressan in \cite[Problem 3, p. 120]{Bressan}. In Bressan's book, the result is 
given for all Lipschitz fluxes. One can take the lower convex envelop of the flux instead of the flux 
in the same formula, thanks to the Oleinik criteria for entropy solutions with general fluxes. 
Moreover, if the flux is strictly convex, then the solution is defined everywhere by formula 
\eqref{b-rarefaction}. Otherwise, it is defined a.e in \cite{Bressan}. 
To be self contained, Proposition \ref{prop:rarefaction} has been added here with a short proof (see also \cite{CC1,CC2}).
\end{remark}

For $n\in\mathbb{N}^*$, let  
$v_i=u_l+(i/n)(u_r-u_l)$ ($i=0,1,\cdots,n$) and let $a_n$ be the sequence of functions such that 
$\forall i$ 
\[ 
a_n(v_i)\ =\ \bar{a}(v_i),
\]
$a_n$ being linear on $[v_i,v_{i+1}]$, so $f_n(u)=f(u_l)+\int_{u_l}^u a_n(v)\,\mathrm{d}v$.  
For proving Proposition \ref{prop:rarefaction}, we need Lemma \ref{lemma1}:
\begin{lemma}\label{lemma1}
For all $(v,\xi) \in [-M,M] \times a([-M,M])$, the sequences $f_n(v)$ and $b_n(\xi)$ converge, 
respectively, to $f(v)$ and $b(\xi)$, and $a_n(v) \to a(v)$ when $n\to\infty $ if $a$ is continuous 
at $v$. Moreover, the sequence $(b_n)$ converges uniformly towards $b$ on any bounded set.
\end{lemma} 
\noindent\textbf{\em Proof.}
\begin{itemize}

\item 
For a given $v$ such that $a$ is continuous at $v$, $v_i=u_l+i (u_r-u_l)/n$ and $v_{i+1}$ are 
chosen such that $v \in [v_i,v_{i+1}]$ (note that $v_i$ depends on $n$). By definition, 
$a_n(v)=n\frac{\bar{a}(v_{i+1})-\bar{a}(v_i)}{u_r-u_l}(v-v_i)+\bar{a}(v_i)$. Since $|v-v_i| \leqslant \frac{u_r-u_l}{n}$ 
and since $a$ is continuous at $v$, then  $\lim\limits_{n \to \infty}a_n(v)=a(v)$. Therefore, 
$\lim\limits_{n \to \infty}f_n(v)=f(v)$ follows at once from the definition of $f_n$.

\item For $n \in \mathbb{N}^* $ and $\xi\in a([-M,M])$, there exists $i$ such that $b(\xi) \in 
[v_i,v_{i+1}]$. Moreover, $b_n$ being linear on $[a_n(v_i),a_n(v_{i+1})]$, $|b_n(\xi)-b(\xi)| 
\leqslant v_{i+1}-v_i=\frac{u_r-u_l}{n}$. In addition, $b$ being continuous, a Dini's Lemma yields 
the uniform convergence of $b_n$ to $b$ as $n\to\infty$.
\cqfd
\end{itemize}

We are now able to prove Proposition \ref{prop:rarefaction}.\\
\textbf{\em Proof.}
By definition of $a_n$, $a_n(u_l)=\bar{a}(u_l)$ and $a_n(u_r)=\bar{a}(u_r)$. Since $f_n\in 
\mathcal{C}^1$, then, from \cite{Bressan, DafermosBook}, the exact entropy solution of 
$u_t+f_n(u)_x=0 $ with the initial condition \eqref{Rm} is
\begin{equation}\label{solution}
u_n(t,x)\,=
\begin{cases}
   u_l & \xi < \bar{a}(u_l), \\
   b_n\left( \xi \right) & \bar{a}(u_l)<\xi< \bar{a}(u_r), \\
   u_r & \bar{a}(u_r)  < \xi. 
\end{cases} 
\end{equation}
Taking the limit $n\to\infty$, thanks to Lemma \ref{lemma1}, $u$ has an explicit continuous formula.  
Thus, $u$ is a weak entropy solution satisfying \eqref{ws} and \eqref{entropy}. Noting that $b(\xi)=u_l$ 
$\forall \xi \in [a^-(u_l),a^+(u_l)]$, $ \bar{a}(u_l)$ in \eqref{solution} can be replaced by $a^-(u_l)$ 
or by $a^+(u_l)$. Similarly, $\bar{a}(u_r)$ can be replaced by $a^-(u_r)$ or  $a^+(u_r)$, thus concluding 
the proof.
\cqfd

Finally, Proposition \ref{prop:whereOlga} can be proved: \\
\textbf{\em Proof.}  
Since $x/t, y/t \in \bar{a}([-M,M])$ and since $b$ is the generalized inverse of $a$, then 
\[ 
\bar{a}\left(u(t,x)\right)=\bar{a}\left(b(x/t)\right)=x/t, \quad \bar{a}\left(u(t,y)\right)=\bar{a}\left(b(y/t)\right)=y/t.
\]
Therefore, the Oleinik inequality is true.
\cqfd


\subsection{Approximate Riemann solver}

In Section \ref{sec4} below, the wave front tracking algorithm is used. Therefore, a suitable 
Oleinik inequality is needed for the approximate solutions. In order to get this inequality, the 
flux is replaced by a suitable piecewise linear approximation.

For this purpose, a piecewise constant approximation of the velocity is used as in \cite{Dafermos72}. A key point is to choose a discrete set of the value of the approximate solution $u^\varepsilon$,  taking into account the discontinuities of the velocity. 
Let $\varepsilon>0$ and let 
\begin{align}
\mathfrak{B}=\{c_0=-M,c_1, \cdots , c_p=M \}
\end{align}  
a subdivision of the 
interval $[-M,M]$ including a too large jump of the velocity. For this purpose the subdivision is chosen such that $c_i<c_{i+1}, \, \forall i$, and
\begin{equation}\label{decomposition}
a^-(c_{i+1})\ -\ a^+(c_i)\ \leqslant\ \quat\,\varepsilon.
\end{equation}
That means that the variation of $a$ is small on $]c_i,c_{i+1}[, \, \forall i$, thus the big jumps of $a$ are located at $c_i$.
\begin{remark}
Due to the jumps of the velocity $ a $ which are not expected on  a uniform grid, the subdivision $ 2 ^ {- n} \mathbb {Z} $ cannot be used as in \cite{Bressan}.
\end{remark}
\begin{remark}\label{largejumps} $\frac{a^+(u)-a^-(u)}{2}= \frac{\llbracket a \rrbracket (u)}{2}$ can be bigger than $\varepsilon/4$, thus,
the condition \eqref{decomposition} cannot be replaced by $\bar{a}(c_{i+1})- \bar{a}(c_i) \leqslant \quat\varepsilon $. Notice that in this case $u$  necessarily belongs to  $\mathfrak{B}$, because the velocity has a big discontinuity at $u$.
\end{remark}
\begin{remark}\label{subdivision}
The condition \eqref{decomposition} is enough to show that if $\varepsilon \to 0$, then $c_{i+1}-c_i \to 0, \, \forall i$, since
\begin{equation}
c_{i+1}-c_i=b \left(a^-(c_{i+1}) \right)-b \left( a^+(c_i)\right) \leqslant \omega[b] \left (a^-(c_{i+1})\ -\ a^+(c_i) \right) \leqslant \omega[b](\quat\,\varepsilon).
\end{equation} 
And $\omega[b]$ is continuous at $0$, thanks to the Heine theorem.
\end{remark}

As in \cite{Dafermos72}, the flux $f$ is approximated by a  continuous and piecewise linear flux. For this purpose, the approximate flux $f_\varepsilon$ is chosen as the  continuous piecewise linear interpolation of $f$ on the subdivision $\mathfrak{B}$, $f_\varepsilon(c_i)=f(c_i) \; \forall i$  and   $f_\varepsilon$  is linear on $[c_i,c_{i+1}]$. Its derivative $a_\varepsilon=f_\varepsilon'$ is piecewise constant.

Now, the approximate Riemann solver is expounded.  Let $u_l,u_r \in \mathfrak{B}$, and let $u^\varepsilon$ be the entropy solution of the Riemann problem $u_t+f_ \varepsilon(u)_x=0$, with the initial data 
\eqref{Rm}.

If $u_l>u_r$, as in Subsection \ref{exactsol}, the solution generates a shock with the Rankine-Hungoniot relation $s=\llbracket f_\varepsilon(u) \rrbracket/\llbracket u \rrbracket=\llbracket f(u) \rrbracket/\llbracket u \rrbracket$ (Notice that $f=f_\varepsilon$ on $\mathfrak{B}$).

If $u_r>u_l$, let $u_l=c_k, \, u_r=c_{k'}$, for a fixed $t>0$, $u^\varepsilon(t,\cdot)$ is non decreasing 
and piecewise constant. Defining 
\begin{equation}\label{si}
s_i\ =\ \frac{f_\varepsilon(c_i)-f_\varepsilon(c_{i-1})}{c_i-c_{i-1}}\ 
=\ \frac{f(c_i)-f(c_{i-1})}{c_i-c_{i-1}}  =a_\varepsilon (c) \mbox{ on } (c_{i-1},c_i).
\end{equation}  
The solution is given by $u^\varepsilon(t,x)=c_i$ for $x/t \in ]s_i,s_{i+1}[$ as in \cite{Dafermos72}. The curves of discontinuity in this case are called contact discontinuities. In fact those curves represent an approximation of a rarefaction wave, so in this paper we call them rarefaction curves. 

Since 
$s_i \leqslant \bar{a}(c_i) \leqslant s_{i+1}$, then \eqref{si} implies
\begin{equation}
u(t,t s_i)\,=\,b(s_i)\,\leqslant\,b(\bar{a}(c_i))\,=\,c_i\,\leqslant\,b(s_{i+1})\,=\,u(t,t s_{i+1}).
\end{equation}
For $i$ with $k<i<k'$, the equation $u(t,\tilde{x}_i)=c_i$ has at least one solution since the exact solution $u$ is non decreasing and continuous. Adding the condition  $\tilde{\xi}_i=\tilde{x}_i/t\in
\bar{a}([-M,M])$, this solution is unique, and  $\tilde{x}_i=\bar{a}(c_i)t $ so $\tilde{\xi}_i
=\bar{a}(c_i)$.

Let $\tilde{\xi}_k,$ $\tilde{\xi}^+_k,$ $\tilde{\xi}_{k'},$ $\tilde{\xi}^-_{k'}$ be defined as
\begin{itemize}
\item $\tilde{\xi}_k=\bar{a}(u_l),$
\item $\tilde{\xi}^+_k=a_\varepsilon^+(u_l),$
\item $\tilde{\xi}_{k'}=\bar{a}(u_r),$ 
\item $\tilde{\xi}^-_{k'}=a^-_\varepsilon(u_r)$.
\end{itemize}
Now, $\tilde{x}_k, \tilde{x}_k^+,\tilde{x}_{k'},\tilde{x}_{k'}^-$ are defined by the relation $x=\xi \, t$.

By construction, the approximate solution at the point $ (t,\tilde{x}_i)$ equals the exact solution 
at the same point, i.e., 
$u^\varepsilon(t,\tilde{x}_i)=u(t,\tilde{x}_i)$ $\forall i=k+1, \ldots , k'-1 $. Since $\tilde{\xi}_i
\in\bar{a}([-M,M])$, then 
\[
\bar{a}(u^\varepsilon(t,\tilde{x}_{i+1}))\ -\ \bar{a}(u^\varepsilon(t,\tilde{x}_i))\ =\ \bar{a}(b(\tilde{x}_{i+1}/t))\ -\ \bar{a}(b(\tilde{x}_i/t))\ =\ 
(\tilde{x}_{i+1}-\tilde{x}_i)/t, \, \forall i=k+1, \ldots , k'-2 .
\]

\begin{center}
\begin{tikzpicture}
\draw[thick,->] (-5,-3) -- (5,-3) node[anchor= west] {$x$};

\draw[thick][color=blue] (0,-3) -- (-4,0) ;
\draw[thick][color=blue] (0,-3) -- (-2,0) ;
\draw[thick][color=blue] (0,-3) -- (0.5,0) ;
\draw[thick][color=blue] (0,-3) -- (3,0) ;

\draw[dashed][color=red] (-0.4,-2.75) -- (-4.1,0) ;
\draw[dashed][color=red] (0.4,-2.7) -- (3.1,0) ;

\draw[dashed] (0,-3) -- (-4.2,-1) ;
\draw[dashed] (0,-3) -- (-3,0) ;
\draw[dashed] (0,-3) -- (-0.75,0) ;
\draw[dashed] (0,-3) -- (1.5,0) ;
\draw[dashed] (0,-3) -- (3.2,-1) ;

\draw (-4.2,-1) node {$\tilde{\xi}_0$};
\draw (-3,0) node {$\tilde{\xi}_1$};
\draw (-0.75,0) node {$\tilde{\xi}_2$};
\draw (1.5,0) node {$\tilde{\xi}_3$};
\draw (3.2,-1) node {$\tilde{\xi}_4$};

\draw (-3.5,-0.6) node[color=red] {$\tilde{\xi}^+_0$};
\draw (3,-0.4) node[color=red] {$\tilde{\xi}^-_4$};

\draw (-4,0) node[above][color=blue] {$s_1$};
\draw (-2,0) node[above][color=blue] {$s_2$};
\draw (0.5,0) node[above][color=blue] {$s_3$};
\draw (3,0) node[above][color=blue] {$s_4$};

\draw (-3,-3) node[above] {\Large $u_l=u_0$};
\draw (2.5,-3) node[above] {\Large $u_r=u_4$};

];
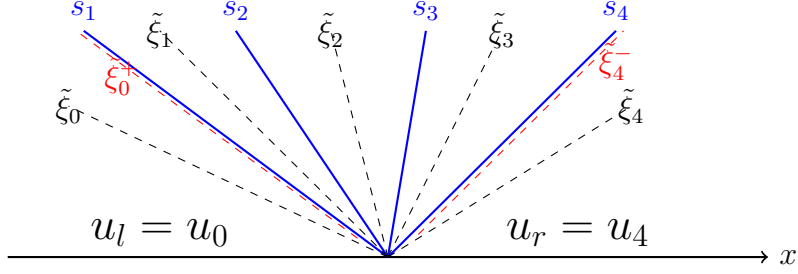
\node at (0,-4) {\textbf{Fig. 2.} An example of a rarefaction wave with $k=0, \, k'=4$.};
\end{tikzpicture}
\end{center}

Summing up for all $i$, since 
\[
a^-(u_{k'})-a_\varepsilon^-(u_{k'}) \leqslant a^-(u_{k'})-a^+(u_{k'-1})\leqslant \varepsilon/4, 
\]
and
\[
a_\varepsilon^+(u_{k})-a^+(u_{k}) \leqslant a^-(u_{k+1})-a^+(u_{k}) \leqslant \varepsilon/4, 
\]
the error is smaller than $\varepsilon/4$ on each boundary term. Thus, the modified Oleinik inequality 
holds also on the whole interval as
\begin{subequations}\label{mOi}
\begin{align}\label{mOia}
\bar{a}(u^\varepsilon(t,\tilde{x}_{k'}))\ -\ \bar{a}(u^\varepsilon(t,\tilde{x}_k))\ &=\ 
(\tilde{x}_{k'}-\tilde{x}_k)/t\,\\ \label{mOib}
a^-(u^\varepsilon(t,\tilde{x}^-_{k'}))\ -\ \bar{a}(u^\varepsilon(t,\tilde{x}_k))\ &\leqslant\ 
(\tilde{x}^-_{k'}-\tilde{x}_k)/t\ +\ \varepsilon/4, \\ \label{mOic}
\bar{a}(u^\varepsilon(t,\tilde{x}_{k'}))\ -\ a^+(u^\varepsilon(t,\tilde{x}^+_k))\ &\leqslant\ 
(\tilde{x}_{k'}-\tilde{x}^+_k)/t\ +\ \varepsilon/4, \\ \label{mOid}
a^-(u^\varepsilon(t,\tilde{x}^-_{k'}))\ -\ a^+(u^\varepsilon(t,\tilde{x}^+_k))\ &\leqslant\ 
(\tilde{x}^-_{k'}-\tilde{x}^+_k)/t\ +\ \varepsilon/2.
\end{align}
\end{subequations}
The value of $u^\varepsilon(t,\tilde{x}^-_{k'})$ is considered on the right of the curve of 
discontinuity, and the value of $u^\varepsilon(t,\tilde{x}^+_{k})$ is considered on its
left (see figure 2), i.e., $u^\varepsilon(t,\tilde{x}^-_{k'})=u_{k'}=u_r, \, 
u^\varepsilon(t,\tilde{x}^+_{k})=u_k=u_l $.
This proves the approximated Oleinik inequality for rarefaction waves.
\cqfd

Here, the rarefaction wave has been approached by a sequence of step functions and satisfies 
the approximate Oleinik inequality \eqref{mOi}. 
Usually,  the Oleinik inequality gives that $a(u)$ is in  $\mathrm{BV}$ and then $u$ in $\mathrm{BV}^\Phi$ 
\cite{CJJ}. Unfortunately, $a(u)$ is not well defined. Moreover, the modified Oleinik inequality \eqref{mOi} 
does not imply  that $\bar{a}(u)$ is $\mathrm{BV}$. 
The next section is devoted to define another velocity  $\chi \cong a(u^\varepsilon)$ which can be controlled 
in  $\mathrm{BV}$ with the wave front tracking algorithm and  the restricted Oleinik inequality \eqref{mOi}.

%

\section{Wave front tracking algorithm}\label{sec4}

This Section deals with the $\mathrm{BV}$ estimate of a velocity $\chi$ defined below. For that purpose, 
a $\mathrm{BV}^+$ estimate is used. With $(x)^+\eqdef\max \{x,0\}$, the $\mathrm{BV}^+$ space is defined 
by $\mathrm{BV}^+\eqdef\{u,\,\mathrm{TV}^+u<\infty\}$, $\mathrm{TV}^+$ being the {\it positive total 
variation}
\begin{equation}
\mathrm{TV}^+\/v\ \eqdef\ \sup_{p\/\in\/\mathcal{P}}\ \sum_{i=2}^n \left(v(x_i)-v(x_{i-1})\right)^+,
\end{equation}
where $\mathcal{P}=\{ \{ x_1,x_2,\cdots,x_n\}, \, x_1<\cdots<x_n, \, 1 \leqslant n \}$ is the set of all 
subdivisions of 
$\mathds{R}$.

The function $u_0$ being bounded, we can assume that $u_0$ has a compact support to prove 
Theorem \ref{maintheorem} (thanks to the finite speed of waves propagation). 

Let $A$ be a positive real number such that $\mathop{\mathrm{supp}}(u_0) \subset [-A,A]$.  Let $ h=A/2m$ and 
$x_i=-A+hi$ with $i=0,1,\cdots,m$. The initial datum $u_0$ is approached by a sequence of step 
functions $(u_{0,m})_m$ taking values in $\mathfrak{B}$ (see remark \ref{subdivision}) and constant 
in $]x_i,x_{i+1}[$ as in \cite[Chapter XIV]{DafermosBook}. 
Consider the initial value problem
\begin{equation}\label{mcl}
u_t\ +\ f_\varepsilon(u)_x\ =\ 0, \qquad u(0,x)\ =\ u_{0,m}.
\end{equation}
The entropy solution of \eqref{mcl}, $u^\varepsilon_m$, is piecewise constant 
\cite{Bressan, DafermosBook}.
The problem \eqref{mcl} requires to solve $m+1$ Riemann problems and the need to study the 
waves interactions.


Note that in the special case of figure 2, if there is a shock on the right of the rarefaction wave, 
that has the values $u_4$ on its left and $u_5 < u_3$ on its right, then the 
distance between the shock and the rarefaction becomes very small. 
There, the total variation of $\bar{a}(u)$ is bigger than 
$\llbracket a \rrbracket (u_4)$, that is problematic because the total variation
 of $\bar{a}(u)$ can not be controlled by the distance between the rarefaction
  wave and the shock wave.
To avoid this problem, the velocity on the part $u=u_4$ is replaced by 
$a^-(u_4)$ instead of $\bar{a}(u_4)$.

In the general case, a new velocity, $\chi^\varepsilon_m$, is introduced. This velocity is defined by removing the jumps of $a(u)$ on the boundaries of the rarefaction wave, if this wave is close to a shock.
Consider $t>0$ a fixed time and $x \in \mathds{R}$. If there is a shock on the left of the point 
$(t,x)$  and a rarefaction on its right, then 
\[ 
\chi^\varepsilon_m(t,x)=a^+(u^\varepsilon_m(t,x)). 
\] 
If there is a shock on the right of the point $(t,x)$ and a rarefaction of its left, then 
\[
\chi^\varepsilon_m(t,x)=a^-(u^\varepsilon_m(t,x)). 
\]
Otherwise 
\[
\chi^\varepsilon_m(t,x)=\bar{a}(u^\varepsilon_m(t,x)).
\]

This definition avoids the problem mentioned above for the special case of figure 2.

Note that in all the three cases, the solution $u^\varepsilon_m$ can be obtained by 
$u^\varepsilon_m(t,x)=b(\chi^\varepsilon_m(t,x))$, which is a key point to take the 
limit in Section \ref{sec5}. Note also that if $f \in \mathcal{C}^1$,  then 
$\chi^\varepsilon_m=a(u^\varepsilon_m)$.

The choice of the inequality \eqref{mOi} depends on the following cases:
\begin{itemize}
\item If two shocks appear on both sides of the rarefaction wave, then the inequality \eqref{mOid} 
is used;
\item If a shock appears only on the left of the rarefaction wave, then the inequality \eqref{mOic} 
is used;
\item If a shock appears only on the right of the rarefaction wave, then the inequality \eqref{mOib} 
is used;
\item Else, the inequality \eqref{mOia} is used.
\end{itemize}
For $t>0$, the positive total variation of $\chi^\varepsilon_m$ over a rarefaction is smaller than 
the size of the rarefaction divided by $t$, plus a small error \eqref{mOi}. For a shock, the positive 
total variation is equal to zero. 

Here, the positive total variation is estimated after wave interactions, as in Bressan's book 
\cite[Chap. 6, Prob. 6]{Bressan}. Let $u_1$, $u_2$ and $u_3$ be the values of $u_m^\varepsilon$ from the 
left to the right. The speed of the left jump is $s_1=\frac{f_\varepsilon(u_2)-f_\varepsilon(u_1)}{u_2-u_1}$ 
and the speed of the right jump is $s_2=\frac{f_\varepsilon(u_2)-f_\varepsilon(u_3)}{u_2-u_3}$. The 
rarefaction wave is replaced by a contact discontinuity since the flux is piecewise affine. 
All the possibilities for wave interactions are listed below.
\begin{enumerate}
\item[(SS)] Shock--shock interaction: $u_3<u_2<u_1$. When two shocks collide they generate a new shock, 
and the positive total variation is always equal to zero.
\item[(RS)] Rarefaction--shock interaction: $u_3 \leqslant u_1<u_2$. the values $u_1$ and $u_2$ are consecutive 
in $I$, which implies that $u_3 \leqslant u_1$. We consider the two cases:
\begin{itemize}
\item $u_1=u_3$. After the interaction, the curves of discontinuity will disappear, and the positive 
total variation will be zero.
\item $u_3<u_1$. After the interaction, a shock will appear (see figure), the rarefaction fan will be 
smaller, and it will lose the curve of discontinuity on the right (between $u_1$ and $u_2$). The value 
of $\chi^\varepsilon_m$ will be changed from $\bar{a}(u^\varepsilon_m(u_1))$ to $a^-(u^\varepsilon_m(u_1))$ 
(see the definition of $\chi^\varepsilon_m$), and also the choice of the points will be changed from 
$\tilde{\xi}_1$ to $\tilde{\xi}^-_1$, which makes the inequality \eqref{mOi} holds after the interaction.
\vspace{7mm}
\begin{center}
\begin{tikzpicture}
\draw[thick][color=blue] (-2,-3) -- (-4,0);
\draw[thick][color=blue] (-2,-3) -- (0,-1.5) ;
\draw[thick][color=blue] (2,-3) -- (0,-1.5) ;
\draw[thick][color=blue] (-1.4,0) -- (0,-1.5) ;
\draw[dashed] (-2.9,-1.5) -- (-3.9,0) ;
\draw[dashed] (-2,-3) -- (-1.8,-1.5) ;
\draw[thick,->] (-5,-3) -- (5,-3) node[anchor= west] {$x$};
\draw (-4,0) node[left][color=blue] {$s_0$};
\draw (-1,-2.25) node[right][color=blue] {$s_1$};
\draw (1,-2.25) node[left][color=blue] {$s_2$};
\draw (-1.75,-2.25) node[left]{$\tilde{\xi}_1$};
\draw (-2.8,-0.5) node[left]{$\tilde{\xi}^-_1$};
\draw (-4,-3) node[above] {\Large $u_0$};
\draw (-1.5,-1.2) node[above] {\Large $u_1$};
\draw (0,-3) node[above] {\Large $u_2$};
\draw (3,-3) node[above] {\Large $u_3$};
];
\node at (0,-4) {\textbf{Fig. 3.} Interaction (RS).};
\end{tikzpicture}
\end{center}

\end{itemize}
\item[(SR)] Shock--rarefaction interaction: $u_2<u_3 \leqslant u_1$. This case can be treated exactly 
like the case Rarefaction--shock.
\item[(RR)] Rarefaction--rarefaction interaction: $u_1<u_2<u_3$. 
Two rarefactions cannot collide (even the points $\tilde{x}_i$). This case is impossible, because the 
convexity of $f_\varepsilon$ implies that $s_1<s_2$.
\end{enumerate} 

\begin{remark}
In the case (RS) the new shock that appears can be very close to the rarefaction, which means that if 
the jump of the velocity $a$ on $u_1$ is big enough, then, the positive total variation of $\bar{a}(
u^\varepsilon_m)$ cannot be controlled by the size of the rarefaction wave. That is the reason of 
using the function $\chi^\varepsilon_m$.
\end{remark}

In summary, the positive total variation of $\chi^\varepsilon_m$ and the number of rarefaction waves 
do not increase. Also, $\mathrm{TV}^+\/\chi^\varepsilon_m(t,\cdot)$ is bounded by summing up  the 
positive variation of all the rarefaction waves, thanks to the modified one-sided Oleinik inequality 
\eqref{mOi}. For each rarefaction wave, this variation is related to the size of the interval 
at time $t>0$ up to a small error $\varepsilon/2$.
The sum of all sizes of the intervals cannot exceed the size of the support of the solution.
Since the number of rarefaction waves is less than $m$, there are only $m$ error terms of size 
less than $\varepsilon/2$. Thence
\begin{equation}\label{TV^+}
\mathrm{TV}^+\/\chi^\varepsilon_m(t,\cdot)\ \leqslant\ L(t)/t\ +\ m\,\varepsilon/2,
\end{equation}
where  $L(t)=2A(t)$ and $\mathop{\mathrm{supp}}u^\varepsilon_m(t,\cdot) \subset [-A(t),A(t)]$.
Recall that $\mathrm{BV}^+ \cap \mathrm{L}^\infty=\mathrm{BV}$ since 
\[
\mathrm{TV}\/ \chi\ \leqslant\ 2\,(\mathrm{TV}^+\chi+\| \chi \|_\infty).
\]
The boundedness of the propagation velocity yields $L(t) \leqslant 2A+2t \|a(u)\|_\infty$. 
Then, taking the constant $C=C(u_0,f)=\max(4A,6\|a(u)\|_\infty) > 0$, which doesn't depend 
on $\varepsilon$ and $m$, gives
\begin{equation}\label{TV}
\mathrm{TV}\/\chi^\varepsilon_m(t,\cdot))\ \leqslant\ C\left(1+1/t\right)\ +\ m\,\varepsilon.
\end{equation}

\section{Compactness and regularity}\label{sec5} 

This Section is devoted to the proof of Theorem \ref{maintheorem}. For this purpose, an uniform 
estimate of the velocity $\chi^\varepsilon_m$ is obtained, 
which gives the compactness of the sequence $(\chi^\varepsilon_m)$. To take the limit in 
\eqref{TV} as $m\to\infty$ and $\varepsilon\to0$, the parameters $m$ can be chosen such that
\begin{equation}\label{epsilon}
\lim_{\varepsilon \to 0}\ m_\varepsilon\,\varepsilon\ =\ 0.  
\end{equation}
We begin with an estimate of the velocity. 

The $\mathrm{BV}$ estimate of the velocity $\chi^\varepsilon_m$ proved in the previous Section 
yields to $\mathrm{Lip}_t\mathrm{L}^1_x$ estimates.
First, let $t$ be a fixed time belonging to $[T_1,T_2]\subset ]0,+\infty[$. We have
\begin{equation}\label{cx}
\int_\mathds{R} \Big| \chi^\varepsilon_m(t,x+h)-\chi^\varepsilon_m(t,x) \Big| \mathrm{d}x\ 
\leqslant\ TV \chi^\varepsilon_m(t,\cdot) |h|\ \leqslant\, \left[ C \left(1+\frac{1}{T_1}\right) + 
m\, \varepsilon  \right] |h|.
\end{equation}
Second, consider two different times $T_1<T_2$. $\chi^\varepsilon_m$ is piecewise constant and 
it has exactly the same curves of discontinuity of $u^\varepsilon_m$. These curves are Lipchitz and 
the speed of any curve cannot exceed $k=\|a(u_0)\|_\infty$. We suppose at first that there is no wave 
interaction between $T_1$ and $T_2$. Then the domain $[T_1,T_2] \times \mathds{R}$ can be divided as 
the following
\begin{center}
\begin{tikzpicture}
\draw[thick] (5,0) -- (-5,0) node[anchor= east] {$\tilde{T}_4=T_2$};
\draw[thick] (5,-3) -- (-5,-3) node[anchor= east] {$\tilde{T}_1=T_1$};
\draw[] (5,-5/3) -- (-5,-5/3) node[anchor= east] {$\tilde{T}_2$};
\draw[] (5,-1/1.2) -- (-5,-1/1.2) node[anchor= east] {$\tilde{T}_3$};
\draw[thick][color=blue] (1.2 ,-3) -- (3,0) ;
\draw[thick][color=blue] (-2,-3) -- (2.5,0) ;
\draw[thick][color=blue] (-9/2 ,-3) -- (0,0) ;
\draw[dashed] (2.5,-1/1.2) -- (2.5,0) ;
\draw[dashed]  (1.5,-5/3) -- (1.5,-1/1.2) ;
\draw[dashed]  (0,-3) -- (0,0) ;
\draw[dashed]  (-2,-3) -- (-2,-5/3) ;
];
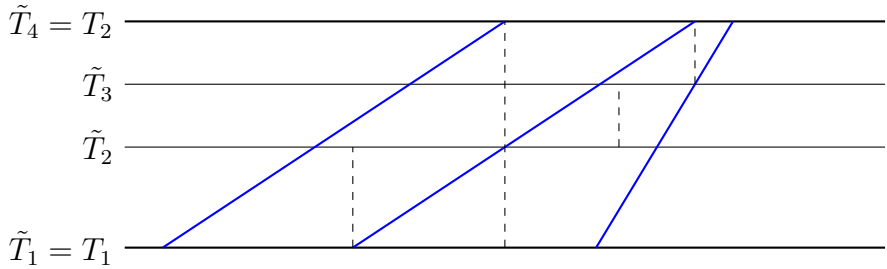
\node at (0,-4) {\textbf{Fig. 4.} Decomposition of the domain, the blue lines are the curves 
of discontinuity.};
\end{tikzpicture}
\end{center}

Using (\ref{cx}) within each small rectangle, gives
\[
\int_\mathbb{R}\Big| \chi^\varepsilon_m(\tilde{T}_j,x)-\chi^\varepsilon_m(\tilde{T}_{j+1},x)
\Big|\mathrm{d}x \leqslant K \left[ C \left(1+\frac{1}{T_1}\right) + m\, \varepsilon \right] 
\Big|\tilde{T}_{j+1}-\tilde{T}_j\Big|.
\]
In general, there are many interactions, so $[T_1,T_2]= \cup_{j=0}^J[t_j,t_{j+1}]$, where $t_0=T_1$,  
$t_J=T_2$ and the points $t_j, \; j=1, \cdots,  J-1$ are the instants of the interactions. Let 
$0<\delta< \half\inf_j (t_{j+1}-t_j)$. The inequality holds true for $t \in [t_j+\delta,t_{j+1}-\delta]$. 
Taking $\delta \longrightarrow 0 $ and using that $\chi^\varepsilon_m(\cdot,x)$ is continuous at 
$t_j$ for almost all $x$, the inequality 
\begin{equation}\label{ct}
\int_\mathbb{R}\Big|\chi^\varepsilon_m(T_1,x)-\chi^\varepsilon_m(T_2,x)\Big|\mathrm{d}x 
\leqslant K \left[ C \left(1+\frac{1}{T_1}\right) + m\, \varepsilon  \right] \Big|T_1-T_2\Big|,
\end{equation}
follows.
Notice that in the estimates \eqref{cx} and \eqref{ct}, the term $m\ \varepsilon$ is bounded by a constant, 
thanks to \eqref{epsilon}. These two inequalities and the uniform boundedness of $\chi^\varepsilon_m$ 
by $k$ are the conditions of the classical compactness Theorem A.8 in \cite{HR}. Hence, the sequence 
$\chi^\varepsilon_m$ converges up to a sub-sequence (if necessary) to some function $\chi$ in $\mathcal{C}
([T_1,T2], \mathrm{L^1_{loc}})$
\begin{equation}
\chi=\lim_{\varepsilon \to 0}\ \chi^\varepsilon_m. 
\end{equation}
Due to the lower semi-continuity of the total variation, we have $\chi \in BV$ and Lipchitz in time 
with value in $\mathrm{L^1_{loc}}$ in space. Thus, for any $0<t$, $\chi$ satisfies
\begin{align}
\mathrm{TV}\,  \chi(t,\cdot) &\leqslant  C \left(1+\frac{1}{t}\right) ,  \label{vBV} \\
\int_\mathbb{R}\Big|\chi(T_1,x)-\chi(T_2,x)\Big|\mathrm{d}x & \leqslant  K  C \left(1+\frac{1}{T_1}\right)   
\left|T_1-T_2\right|. \label{vLip}
\end{align}
Using that $u^\varepsilon_m=b(\chi^\varepsilon_m)$, which also provides the compactness of the 
sequence of approximate solutions. Taking the limit $m\to\infty$ in \eqref{ws} and \eqref{entropy}, 
gives that 
\begin{equation} \label{u=bg}
u=b(\chi)
\end{equation} is an entropy solution of \eqref{cl}.  
The main Theorem of \cite{GM} (see also \cite{P05,P07}) ensures that the initial datum  is recovered. 
Then, $u$ is the unique Kruzkov entropy solution with the initial datum $u_0$.
\begin{remark}
Equality \eqref{u=bg} means that $\chi=a(u)$ for a smooth velocity. Here, it is not necessarily true  
almost everywhere since the velocity $a(u)$ can be discontinuous where $u$ is constant.
\end{remark}
Now, the $\mathrm{BV}^\Phi$ regularity of the entropy solution is proven. Let us first check that 
$\Phi$ is positive. 
The function $b$ is not constant on  the whole interval $[-M,M]$. Then, $\omega(h)>0$ for $h > 0$. 
Thanks to Heine theorem, $\omega$ is continuous at $0$, ensuring that $\phi(y)>0$ for $y>0$.  
$\Phi$ the convex envelope of $\phi$ is then also strictly positive \cite{CJLO}.

The $\mathrm{BV}^\Phi$ regularity of $u$ is a direct consequence of the $\mathrm{BV}$ regularity of 
$\chi$ and the definition of $\Phi$, which yield with \eqref{u=bg}, \eqref{modulus}, and remark 
\ref{leftinverse} to the following inequality for almost all $t_1,t_2,x,y$

\begin{align}\label{phi}
\Phi(| u(t_1,x)-u(t_2,y) | )\ &=\ \Phi(| b(\chi(t_1,x))-b(\chi(t_2,y)) | ) \nonumber  \\
&\leqslant\ \phi(|b(\chi(t_1,x))-b(\chi(t_2,y))|)     \nonumber  \\
&\leqslant\ \phi \Big( \omega \big( |\chi(t_1,x) -\chi(t_2,y)| \big)  \Big)   \nonumber  \\ 
&\leqslant\ \big| \chi(t_1,x) -\chi(t_2,y)  \big|.
\end{align}
The $\mathrm{BV}$ regularity of $\chi$ \eqref{vBV} and the inequality \eqref{phi} show that 
$u \in BV^\Phi$. The $\mathrm{Lip}_t \mathrm{L}^1_x$ regularity of $\chi$ \eqref{vLip} and 
inequality \eqref{phi} again implies that 
\begin{equation}
\int_\mathbb{R} \Phi \Big( \big|u(t_1,x)-u(t_2,x) \big| \Big) \mathrm{d}x \leqslant KC(1+\tau^{-1})|t_1-t_2|.
\end{equation}
This is an estimate in the Orlicz space $\mathrm{L}^\Phi$, i.e.,   
$\mathrm{L}^\Phi(\mathds{R})$ denotes the set of measurable functions $f$ such that 
$\int_\mathds{R}\Phi(|f(x)|) \/\mathrm{d}x < \infty$ \cite{Orlicz}.    

If $\Phi(u)=|u|^p$ then the $\mathrm{Lip}^\frac{1}{p}\left(]\tau,+\infty[, 
\mathrm{L}^p_\mathrm{loc}\right)$ estimate in \cite{BGJ} is recovered. In general, 
Jensen inequality gives $u\in \mathcal{C}_w(]\tau,+\infty[, \mathrm{L^1_{loc}})$.\\

\textbf{Proof of the inequality \eqref{mt+}.}

The proof of Theorem \ref{maintheorem} is done for  $u_0$  compactly supported. In the general case, 
to prove the inequality \eqref{mt+}, the segment $[\alpha,\beta]$ is divided by a subdivision 
$\mathfrak{S}=\{\alpha_0=\alpha, \alpha_1, \cdots, \alpha_l=\beta \}$  in order to 
separate rarefaction waves and shocks. Inequality \eqref{mOi} implies that 
\begin{equation}
a^-(u(t,\alpha_{i+1}))-a^+(u(t,\alpha_i))\ \leqslant\, \left( \alpha_{i+1}-\alpha_i \right)t^{-1}
\ +\ \varepsilon/4 \ . 
\end{equation}
Replacing $\chi(t_1,x) -\chi(t_2,y)$ by $ a^-(u(t_1,x))-a^+(u(t_2,y)) $ in the inequality \eqref{phi}, 
the result follows by restarting a similar proof, summing up for all $i$ and taking the limit $m\to +\infty$.






\end{document}